\documentclass{amsproc}
\usepackage{amscd, eucal}

\begin{document}

\newsymbol \twoheadrightarrow 1310
\newsymbol \rightarrowtail 131A

\newcommand{\B}{{\bf B}}
\newcommand{\Br}{{\rm Br}}
\newcommand{\C}{{\mathbb C}}
\newcommand{\config}{{\rm Config}}
\newcommand{\End}{{\rm End}}
\newcommand{\Hom}{{\rm Hom}}
\newcommand{\M}{{\overline{\mathcal M}}}
\newcommand{\Oh}{{\mathcal O}}
\newcommand{\R}{{\mathbb R}}
\newcommand{\Q}{{\mathbb Q}}
\newcommand{\T}{{\mathbb T}}
\newcommand{\PGl}{{\rm PGl}}
\newcommand{\Z}{{\mathbb Z}}

\title {Braids, trees, and operads}\bigskip
\author{Jack Morava}
\address{Department of Mathematics, Johns Hopkins University,
Baltimore, Maryland 21218}
\email{jack@math.jhu.edu}
\thanks{The author was supported in part by the NSF}
\date {5 June 2001}
\begin{abstract} The space of unordered configurations of distinct points
in the plane is aspherical, with Artin's braid group as its fundamental 
group. Remarkably enough, the space of ordered configurations of 
distinct points on the real projective line, modulo projective equivalence,
has a natural compactification (as a space of equivalence classes of trees)
which is also (by a theorem of Davis, Januszkiewicz, and Scott) aspherical.
The classical braid groups are ubiquitous in modern mathematics, with 
applications from the theory of operads to the study of the Galois group
of the rationals. The fundamental groups of these new configuration 
spaces are not braid groups, but they have many similar formal properties.
This talk [at the Gdansk conference on algebraic topology 05-06-01] is an 
introduction to their study. \end{abstract}

\maketitle

\section{The bubbletree operad and quantum cohomology} \bigskip

\noindent
{\bf 1.0} Work on {\bf conformal field theories} leads physicists to an
interest in configuration spaces $$\config^{n+1} \C P_1 \sim \config^n 
\C \;.$$ of points on the complex projective line. They are most interested 
in the {\bf quotients} of these spaces by the action of $\PGl_2(\C)$. 
The points are noncoincident, so both the spaces and the group are
noncompact, and taking the quotient is tricky: it leads naturally
to a compactification $$\M_{0,n}(\C) \sim \config^n(\C P_1)/\PGl_2(\C) \;.$$
The physicists discovered a {\bf repulsive potential} among these
points: pushing two together creates a bubble onto which they
escape [13]. \bigskip

\noindent
{\bf 1.1} Thus $\M_{0,n}(\C)$ is the moduli space of marked genus zero 
{\bf stable} algebraic curves (which have (at worst) double points, 
and at least three marked points on each irreducible component). 
\bigskip

\noindent
{\bf Example:} $\M_{0,4}(\C) \cong \C P_1$ via the classical cross-ratio. 
Note, $\M_{0,3}(\C)$ is a point: a configuration of three points on 
$\C P_1$ is {\bf rigid}. \bigskip

\noindent
These spaces are very nice in some ways: they are compact
manifolds, with cohomology concentrated in even dimension, and no
torsion [8]. \bigskip

\begin{center}
{\bf Operads, by example:}
\end{center} \bigskip

\noindent
{\bf 1.2} An operad $\Oh_* = \{ \Oh_k, k \geq 1 \}$ is a collection of
spaces together with some {\bf composition} maps
$$\Oh_n \times \Oh_{i_1} \times \dots \times \Oh_{i_n} \to \Oh_i$$
(where $i = \sum i_k$) satisfying some axioms \dots \medskip

{\bf ex. i)} $\M_{0,*+1}(\C)$ \medskip

{\bf ex. ii)} $\End_n(X) = {\rm Maps}(X^n,X)$ \medskip

\noindent
is the {\bf endomorphism operad} of an object $X$ in a monoidal category.
Composition sends a map from $X^n$ to $X$ and a tuple of maps from 
$X^{i_*}$ to $X$ to a composition $$X^i = X^{i_1} \times \dots X^{i_n} \to 
X^n \to X \;.$$ A  morphism $\Oh_* \to \End_*(X)$ of operads makes $X$
into an $\Oh_*$-{\bf algebra}. \medskip

{\bf ex. iii)} $\Br_n$ = Artin's braid group on $n$ strings
\medskip

\noindent
defines the {\bf braid operad}, with {\bf cabling} 
$$\Br_n \times \Br_{i_1} \times \dots \Br_{i_n} \to \Br_i$$
as composition. \bigskip

\noindent
Monoidal functors preserve operads; hence the homology of an operad 
(in spaces) is an operad in graded modules [9].
\bigskip

\noindent
{\bf 1.3 Theorem} (WDVV, Kontsevich [5]): The (rational) homology
of a smooth projective algebraic variety $V$ is an
$H_*(\M_{0,*+1}(\C))$-operad algebra. \medskip

\noindent
[This led to the solution of the 19th-century enumerative
geometry problem of classification of lines with specified
incidence in $\C P_2$.] \bigskip

\noindent
The {\bf construction} of this algebra structure uses the 
Gromov-Witten invariants: There are (compact) moduli spaces 
$GW_k(V)$ of holomorphic maps from genus zero stable curves with 
$k$ marked points, to $V$. These spaces have many components, 
indexed by degree $[h:C \to V] \in H_2(V,\Z)$. There is also 
an {\bf evaluation} map $$GW_k(V) \to \M_{0,k}(\C) \times V^k$$ 
which defines a cycle $$GW_k \in H_*(\M_{0,k}(\C)) \otimes 
H_*(V)^{\otimes k} \;.$$ [Actually the coefficients lie in the 
Novikov ring $\Lambda = \Q[H_2(V,\Z)]$, but this will be supressed.] 
Using Poincar\'e duality, we can rewrite $GW_{k+1}$ as an element of 
$$\Hom(H_*(\M_{0,k+1}(\C)),\Hom(H_*(V)^{\otimes k},H_*(V)))$$
which then defines a morphism $$H_*(\M_{0,k+1}(\C)) \to \End_k(H_*(V))$$
of operads, QED. \bigskip

\noindent
In particular, the point $\M_{0,3}(\C)$ defines a {\bf quantum 
multiplication} $$H_*(V,\Lambda) \otimes_{\Lambda} H_*(V,\Lambda) \to
H_*(V,\Lambda)$$ which is usually not standard \dots \bigskip

\section{Devadoss's mosaic operad and Fukaya's Lagrangian cohomology}

\noindent
{\bf 2.0} The moduli space [2] $$\config^n(\R P_1)/\PGl_2(R) \sim 
\M_{0,n}(\R)$$ of configurations of points on the circle can be pictured 
as a space of {\bf trees} or {\bf mosaics} of hyperbolic polygons. 
The points have an intrinsic cyclic order, and $\{ \M_{0,*}(\R) \}$ 
is naturally a {\bf cyclic} operad [6]. I am indebted to Sasha Voronov,
for pointing out that the analogous complex operad is also cyclic!
\bigskip

\noindent
{\bf 2.1 Fukaya} considers a compact symplectic manifold $(M,\omega)$
together with an oriented Lagrangian submanifold $L$ (i.e. of
half the dimension of $M$, such that $\omega|_L = 0$; some subtle
issues involving the Stiefel-Whitney class $w_2(L)$ will be
ignored.) I {\bf conjecture} that the following is a theorem; something 
slightly weaker (cf. below) has been proved by Fukaya and his school [4]:
\bigskip

\noindent
For a generic almost-complex structure compatible with $\omega$, there are
compact oriented moduli spaces $FO_k$ of pseudo-holomorphic hyperbolic 
polygons $$(P,\partial P) \to (M,L)$$ together with evaluation maps
$$FO_k \to \M_{0,k}(\R) \times L^k$$ which define an action of 
$H_*(\M_{0,*+1}(\R))$ on $H_*(L,\Lambda)$ (where now $\Lambda = 
\Q[H_2(M,\Z)]$). \bigskip

\noindent {\bf 2.2} The (co)homology of these spaces is {\bf not} known
(but cf. [14]). However, we can draw some beautiful pictures [2,3]. \bigskip

\noindent
Grothendieck, in his {\it Esquisse}, says $\M_{0,5}(\C)$ is `un petit joyaux'.
Its real points $\M_{0,5}(\R)$ map to $\M_{0,4}(\R) \times
\M_{0,4}(\R) = T^2$ by selecting two distinct subsets of four
points. To get the full space, we need to blow up (ie, add
crosscaps) at the three configurations defined by triple coincidences,
resulting in  $T^2 \# 3 \R P^2$. \bigskip

\noindent
There is a more symmetric picture, defined by blowing up four
points on $\R P^2$. Both pictures are {\bf tesselated} by
pentagons, though this is easier to see in the second picture. 
This is a regular polytope with twelve pentagonal faces: it is
the dodecahedron's `evil twin'. \bigskip

\noindent
{\bf 2.3} In general, there is a surjective map
$$\Sigma_k \times_{D_k} K_{k-3} \to \M_{0,k}(\R)$$
(where $D_k$ is the dihedral group of order $2k$) which is $2^n$ to 
$1$ on codimension $n$ faces: in general, these moduli spaces are 
tesselated by Stasheff {\bf associahedra}.\bigskip

\noindent
There is a commutative diagram
$$
\begin{CD} \config^*(\R) &@>>>& \config^*(\C) \\@VVV &&@VVV\\ \M_{0,*+1}(\R)
&@>>>& \M_{0,*+1}(\C) \end{CD}
$$
The space in the upper right corner is homotopy-equivalent to the
little disks operad, and the space in the upper left corner is
the classical $A_{\infty}$ operad $\{\Sigma_* \times K_{*-1} \}$
(made {\bf permutative}, ie endowed with an action of the
symmetric group. The left vertical map defines the tesselation;
thus the mosaic operad is a kind of (quasicommutative) quotient
of the $A_{\infty}$ operad. Fukaya shows that the $A_{\infty}$ 
operad acts on Floer cohomology, but I believe that action passes
through this quotient. The diagram above is a fiber product of spaces, 
but it is not quite a fiber product of operads.\bigskip 

\noindent
{\bf 2.4 Theorem} of Davis, Januszkiewicz, and Scott [1]: this tesselation
defines a piecewise negatively curved metric on $\M_{0,*}(\R)$;
these spaces are therefore $K(\pi,1)$'s! \bigskip

\noindent
{\bf Remark:} Devadoss [3] has recently shown (using similar methods) that
the Fulton-MacPherson compactification of $\config^*(\R P_1)$ is also 
aspherical! \bigskip 

\section{Operads in groups (and groupoids)} \bigskip

\noindent
{\bf 3.0 Observation:} The fundamental group of an operad is an operad in 
groups \dots provided you're careful about basepoints. \bigskip

\noindent
{\bf Example} $\{1,\dots,n \in \C \}$ (with that order) defines a
basepoint $* \in \config^n(\C)$; but the natural action of
$\Sigma_n$ moves it around (by changing the order). \bigskip

\noindent
Recall that a space has a fundamental {\bf groupoid}, with
respect to a system of basepoints: it is a category, with the
points as objects, and homotopy classes of paths between them as
morphisms. \bigskip

\noindent
Note also that a surjective homomorphism $\phi: G \to H$ of
groups defines a groupoid $[H/G]$ with $H$ as set of objects, and
$${\rm mor}(h_0,h_1) = \{ g \in G \;|\: \phi(g)h_0 = h_1 \}$$
as morphisms. With this notation,
$$\pi(\config^n(\C) \; {\rm rel} \; \Sigma_n(*)) \cong
[\Sigma_n/\Br_n]$$
where $\Br_n \to \Sigma_n$ is the standard homomorphism. Thus the 
fundamental groupoid of the little disks operad is the braid
operad. I owe thanks to Dan Christensen and J.P. Meyer, for explaining
this to me. \bigskip 

\noindent
{\bf 3.1 Definition:} The braid {\bf category} $\B$ has integers $n$ as
objects, with $\Br_n$ as its endomorphisms. [There are no
morphisms between distinct integers.] This is a (universal)
braided monoidal category, with tensor product $\B \times \B \to \B$ 
defined by juxtaposition $(n,m) \mapsto n+m$. \bigskip

\noindent
A standard construction [7] defines a functor
$$[\Sigma_n/\Br_n] \to {\rm Func}(\B^n,\B)$$
which makes the category $\B$ an {\bf algebra} over the braid operad. More
generally, any braided monoidal category is an algebra over the
braid operad. \bigskip

\noindent
The operad $\M_{0,*+1}(\R)$ defines a similar category: there is
an exact sequence
$$\pi_1(\M_{0,*+1}(\R)) \rightarrowtail \pi_1(\M_{0,*+1}(\R)_{h\Sigma_*}) 
\twoheadrightarrow \Sigma_*$$
in which the fundamental group of the homotopy quotient plays the
role of the braid group. There is a similar tensor category {\bf
D}, which is a kind of universal example of an algebra over the
associated operad in groupoids. \bigskip

\noindent
{\bf 3.2} Here are some questions and speculations: \bigskip

\noindent
i) do these fundamental groups act in some natural way on
Fukaya's cohomology? \medskip

\noindent
ii) does {\bf D} have an interpretation in terms of cyclic 
operads with a trace in the sense of Markl [11]? \medskip

\noindent
iii) are these fundamental groups Galois groups for
solutions of Calogero-Moser systems [of points moving on the
line [12]] analogous to the role played by the braid groups in the
Knizhnik-Zamolodchikov equations? \medskip

\noindent
iv) The Grothendieck-Teichm\"uller group [10] is in some sense the
automorphisms of the braid operad. Do automorphisms of the
operad in groupoids defined by the moduli spaces $\M_{0,*}(\R)$
have any similar properties? \medskip

\noindent
iv) Does the rank of $H_1(\M_{0,k+1}(\R))$ equal $\binom{n}{3}$?

\bibliographystyle{amsplain}

\end{document}